# Flexible Demand Resource Pricing Scheme: A Stochastic Benefit-Sharing Approach

Zhaohao Ding    Feng Zhu    Yajing Wang    Ying Lu    Lizi Zhang
Member, IEEE   Student Member, IEEE   Student Member, IEEE   Student Member, IEEE   Member, IEEE
North China Electric Power University
Beijing, 102206, China
[zhaohao.ding, feng.zhu, yajing.wang, yinglu] @ncepu.edu.cn

*Abstract*–With the rapidly increased penetration of renewable generations, incentive-based demand side management (DSM) shows great value on alleviating the uncertainty and providing flexibility for microgrid. However, how to price those demand resources becomes one of the most significant challenges for promoting incentive-based DSM under microgrid environments. In this paper, a flexible demand resource pricing scheme is proposed. Instead of using the utility function of end users like most existing literatures, the economic benefit of flexible demand resources is evaluated by the operation performance enhancement of microgrid and correspondingly the resource is priced based on a benefit sharing approach. An iteration-based chance-constrained method is established to calculate the economic benefit and shared compensation for demand resource providers. Meanwhile, the financial risks for the microgrid operator due to uncertain factors are mitigated by the chance-constrained criterion. The proposed scheme is examined by an experimental microgrid to illustrate its effectiveness.

I. NOMENCLATURE

*A. Indices and Sets*

| | |
|---|---|
| $L, l$ | Number and index of conventional units. |
| $W, w$ | Number and index of wind farms. |
| $S, s$ | Number and index of solar stations. |
| $T, t$ | Number and index of time slots. |
| $C, c$ | Number and index of types of flexible demand resources. |
| $K_c, k_c$ | Number and index of each type of flexible demand resources. |

*B. System Parameters and Functions*

| | |
|---|---|
| $P_l^{con,min}$ | Minimum power output of conventional unit $l$. |
| $P_l^{con,max}$ | Maximum power output of conventional unit $l$. |
| $RD_l$ | Ramp-down rate of conventional unit $l$. |
| $RU_l$ | Ramp-up rate of conventional unit $l$. |
| $MU_l$ | Minimum-up time of conventional unit $l$. |
| $MD_l$ | Minimum-down time of conventional unit $l$. |
| $SU_l$ | Start-up cost for conventional unit $l$. |
| $SD_l$ | Shut-down cost for conventional unit $l$. |
| $a_l$ | Cost for conventional unit $l$ if it is on. |
| $b_l$ | Fuel cost constant of conventional unit $l$ in time $t$. |
| $c_l$ | Fuel cost constant of conventional unit $l$ in time $t$. |
| $R_t$ | Microgrid operation reserve requirement. |
| $P_{w,t}^{wind,\xi}$ | A random parameter to indicate the wind power output of wind farm $w$ in time slot $t$. |
| $P_{s,t}^{solar,\xi}$ | A random parameter to indicate the solar power output of solar station $s$ in time slot $t$. |
| $f_{k_c,t}$ | Original flexible demand resource $k_c$ in time slot $t$. |
| $RD_{k_c}$ | Ramp-down rate of flexible demand resource $k_c$. |
| $RU_{k_c}$ | Ramp-up rate of flexible demand resource $k_c$. |
| $D_{k_c}^{min}$ | Minimum power of flexible demand resource $k_c$. |
| $D_{k_c}^{max}$ | Maximum power of flexible demand resource $k_c$. |
| $T_{k_c}^{min}$ | Minimum continuous-off time of flexible demand resource $k_c$. |
| $T_{k_c}^{max}$ | Maximum continuous-on time of flexible demand resource $k_c$. |
| $T_{k_c}^{op}$ | Availability period of flexible demand resource $k_c$. |
| $N$ | Total number of scenarios in set $S$. |
| $\gamma_t$ | Penalty cost per unit of curtailed wind and solar power in time slot $t$. |
| $\beta$ | Confidence level of the stochastic problem. |

*C. Decision variables*

| | |
|---|---|
| $P_{t,l}^{con}$ | Power output of conventional unit $l$ in time slot $t$. |
| $o_{t,l}$ | Binary variable to indicate if conventional generator $l$ is on in time slot $t$. |
| $u_{t,l}$ | Binary variable to indicate if conventional generator $l$ is started up in time slot $t$. |
| $v_{t,l}$ | Binary variable to indicate if conventional generator $l$ is shut down in time slot $t$. |
| $P_t^{wind}$ | Wind Power utilized in time slot $t$. |
| $P_t^{solar}$ | Solar Power utilized in time slot $t$. |
| $\alpha_{k_c,t}$ | Binary variable to indicate if flexible demand resource of customer $k_c$ is on in time slot $t$. |
| $D_{k_c,t}$ | Power output of flexible demand resource $k_c$ in time slot $t$. |
| $E_t^\xi$ | Amount of curtailed renewable power in time slot $t$. |
| $\pi_{k_c}$ | Compensation for flexible demand resource $k_c$ for per unit contribution in DSM. |
| $\Delta D_{k_c}$ | Change of controlled flexible demand resource |

$k_c$.

## II. INTRODUCTION

The rapid development of renewable energy resource (RES) results in significant intermittence and uncertainties in the operation of power system, particularly for smaller system such as microgrid[1]. Assisted with optimization algorithms, demand side resources and conventional units jointly provide effective means to reduce the impact of the uncertainties of various RES on the operation of microgrid. Therefore, DSM is playing a critical role on providing flexibility for microgrid operation.

According to the United States Department of Energy, demand side management is typically motivated either by pricing signals or incentive payments [2]. Correspondingly, DSM programs can be divided into two basic categories, i.e., price-based DSM, and incentive-based DSM. In price-based DSM, the end-use customers would adjust its demand based on the time-varying price, which has been discussed in many literatures [3-8]. In incentive-based DSM programs, end-use customers act as controllable loads dispatched by external signals. Meanwhile, those demand resource providers would receive monetary compensation as rewards. Numbers of researches have been conducted to investigate the application of incentive-based DSM in the microgrid [9-14]. Contrasted with price-based DSM, incentive-based DSM has greater and faster responsive speed to solve problems of uncertainties[10], which make it playing a more and more important role in the power system operation [15].

One of the most significant challenges for incentive-based DSM is quantitatively determining the compensation for demand resource providers. Most of the existing literature adopts the utility function or comfort level function to represent the willingness of demand resource providers on participating in the incentive-based DSM program. For example, thermal comfort constraints of customers are adopted in air conditioner system as DSM constants. Reference [9] applies the thermal function between room temperature and energy consumption of air conditioners to express the controllable capacity of it. Similarly, for air conditioners, reference [16] applies the length of forced closing time to decide the impact on the comfort of each air conditioning user. Reference [17] applies the Fanger index as a realistic measure for thermal comfort with the ASHRAE 55 standard to evaluate the range for thermal comfort. Reference [18] applies the contracted allowed minimum comfort violation limits of the demand resource providers to limit the comfort level. Moreover, the utility function of demand resource providers is applied to estimate the participation rate of end-use customers. Reference [10] applies a parameter indicating how much controllable loads can be cut off to express the utility. Reference [13] applies a grading scheme with a parameter to classify the demand resource providers according to their desire/readiness to participate in the DSM. Reference [19] investigates different levels of customers' participation rate in DSM and effects on microgrid operation costs. As indicated by the aforementioned literature, most of those researches use the utility function or comfort function to determine the compensation. However, it can be extremely difficult to precisely derive the utility function or comfort function in practical cases.

Besides, there are numbers of uncertain factors in the operation of microgrid. Reference [20] presents a developed two-stage model to consider the wind and PV powers uncertainties. Reference [8] incorporates demand response and energy storage system into the power system to reduce the influence of uncertainties of wind energy. Reference [21] evaluates the operation and significance of energy storage considering the uncertainties of weather conditions. Reference [22] proposes a robust optimization based approach for optimal multi-microgrid operation in a residential scenario considering renewable energy uncertainties. Most of those works focus on the uncertainties in the optimal operation of microgrid. However, determining the compensation for demand resources needs to calculate the economic benefit contributed in the entire pricing period, which involved a large number of uncertain parameters. Correspondingly, different realization of operation parameters would result in different amount of economic benefit, affecting the calculation of compensation.

In this paper, a flexible demand resource pricing scheme based on benefit-sharing is proposed to determine the compensation for incentive-based DSM programs without knowing the utility function of demand resource providers. In the proposed scheme, the economic benefit induced by flexible demand resource is shared by the microgrid operator and flexible demand resource providers and the compensation is determined correspondingly. An iteration-based chance-constrained method is proposed to address the impact of uncertain parameters on economic benefit and mitigate the risk of the microgrid operator for overcompensation

The main contributions of this paper can be summarized as follows:

- A flexible demand resource pricing scheme is proposed in this paper. The compensation for the demand resource providers can be determined without knowing the utility or comfort function of end-users.
- An iteration-based chance-constrained method is established to evaluate the economic benefit of flexible demand resources while mitigating the risk of financial losses for the microgrid operator under stochastic environment.
- Different types of flexible demand resources are examined by the proposed pricing scheme to demonstrate its effectiveness and robustness.

The rest of this paper is organized as follows. The framework model of proposed flexible demand resource pricing scheme is described in Section III. The mathematical model is formulated in Section IV. The experiment case studies and results analysis are provided in Section V, followed by conclusions in Section VI.

## III. MODEL DESCRIPTION

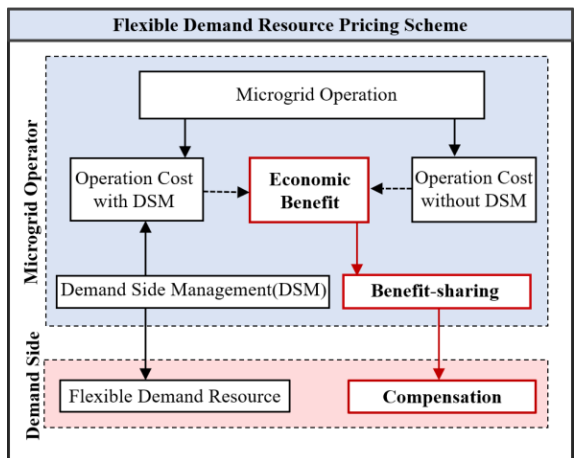

Fig. 1. Structure of flexible demand resource pricing scheme

As mentioned above, a pricing scheme for incentive-based DSM programs under the microgrid environment is proposed, as shown in Fig. 1. The economic benefit of flexible demand resource is measured by the operation cost deviation between adopting or not adopting the incentive-based DSM program. The microgrid operator will determine the compensation (i.e. price) for demand resource providers based on the benefit-sharing scheme.

### A. Demand side management model

In this work, customers have voluntary options to determine whether to participate in the DSM programs. Once the customers sign up the incentive-based DSM programs, flexible loads could be controlled by the microgrid operator. In return, economic benefit saved in the process of DSM program will be distributed to them. In this paper, we try to derive a uniform per unit compensation for certain types of resources to attract more participation.

To motivate more end-users to participate in DSM programs, the compensation distributed to customers will be determined by the specific contribution on economic performance enhancement of microgrid operation. Meanwhile, to avoid potential retail revenue loss for the microgrid operator due to the consumption reduction, all the flexible demand resources considered in this paper follow the "energy-neutral" constraint. In other words, the effect of DSM would only result in load shifting with no reduction of consumption. In this way, the proposed pricing scheme can achieve the Pareto improvement which motivates both the microgrid operator and demand resource providers.

### B. Microgrid model

Generally speaking, microgrid can be categorized into two types depending on whether connecting to the utility grid or not[23]. A microgrid disconnected to the utility grid operates as an islanded microgrid, while connected to the bulk power grid is called a grid- microgrid [24]. In this paper, the islanded microgrid model is adopted to evaluate the economic benefit of flexible demand resources. However, it can be easily extended to grid-tied microgrid by adding the electricity price of utility grid to the proposed model. As shown in Fig. 2, the supply side of microgrid includes conventional generators and renewable resources. The demand side is composed of inflexible demand resources and flexible demand resources.

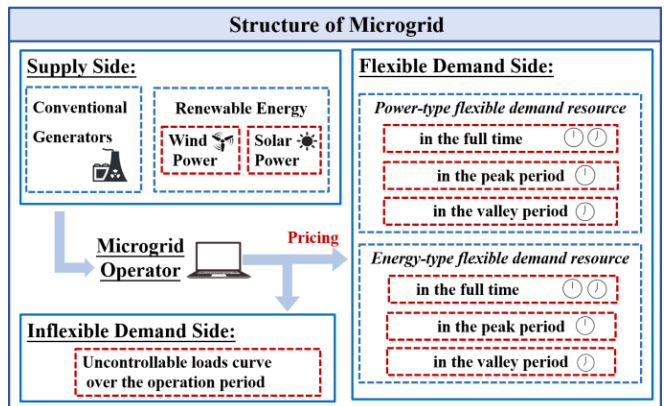

Fig. 2. Structure of the microgrid model

### C. Flexible demand resource model

The flexible demand resource in incentive-based DSM programs can be divided into two categories, power-type flexible demand resource and energy-type flexible demand resource.

As shown in Fig. 3, the power-type demand resource is controllable load with specific power curve in the operation period. Consequently, the microgrid operator could only control its start-up time. That is to say the demand can be shifted however the power curve pattern would stay the same. One of most typical examples for power-type flexible demand resource is industrial processes which follow the fixed power curve once started, such as some processes in the metal processing industry [25].

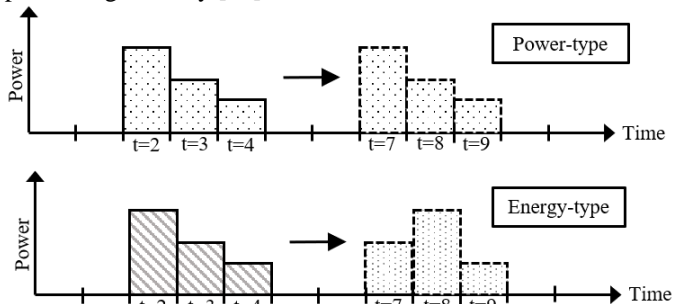

Fig. 3. Examples of power-type and energy-type demand resource

In contrast to that, the energy-type demand resource does not necessarily follow a specific power curve. Instead, the main constraint for it is that certain amount of energy will be consumed in the operation time span, i.e. "energy-neutral". In some cases, it also follows constraints such as ramping up/down and minimum on/off time. A typical example is delay-tolerant cloud computing workloads for data center industry [26]. As long as the workloads can be solved within a time limit, there is no need to follow a fixed power consumption curve. However, the total energy consumed for a certain cloud computing task is generally a constant.

Different availability periods are also considered to precisely describe the characteristics of demand resources. According to classification standards issued by National Development and Reform Commission [27], the operation time span is divided into two periods, the peak period (10:00-20:00) and the valley period (0:00-9:00, 21:00-23:00). Consequently, in this paper, incentive-based DSM programs are divided into three types: 1) power-type/energy-type flexible demand resource available for the entire operation time span; 2) power-type/energy-type flexible demand resource available in the peak period; 3) power-type/energy-type flexible demand resource available in the valley period. The detailed formulation is described in Section IV

## IV. PROBLEM FORMULATION

In this section, the mathematical formulation of flexible demand resource pricing scheme is presented. As indicated by Fig. 1, the compensation for demand resource providers is determined based on the economic benefit induced by flexible demand resources. To calculate the economic benefit based on the microgrid operation model, an iteration-based chance-constrained method is proposed in this paper. The detailed mathematical formulations are provided in the following.

### A. Microgrid operation model

The microgrid operation model is formulated as a classic unit commitment problem with uncertain renewable generations [28], the operation cost consists of the cost of conventional generators and penalty cost of curtailed wind and solar power, which can be expressed as (1) and (2), where the cost of conventional generators is composed of start-up cost, shut-down cost and operation cost. The total power output of conventional generators is described as (3):

$$Cost = \sum_{t=1}^{T}\sum_{l=1}^{L} Cost_{t,l}^{con} + Penalty^{RES} \qquad (1)$$

$$\begin{aligned} Cost_{t,l}^{con} &= SU_l u_{t,l} + SD_l v_{t,l} + a_l o_{t,l} \\ &+ b_l \times P_{t,l}^{con} + c_l \times (P_{t,l}^{con})^2 \quad (\forall l \in L, \forall t \in T) \end{aligned} \qquad (2)$$

$$P_t^{con} = \sum_{l=1}^{L} P_{t,l}^{con} \quad (\forall t \in T) \qquad (3)$$

The problem of unit commitment must meet the constraints listed as follows [29, 30]. The generation capacity constraint of conventional generators is defined in (4). Constraints (5) and (6) represent the minimum-up/down time while the start-up and shut-down constraints are modeled in (7) and (8). The ramping up/down constraints of units are described in (9) and (10). The power balance and system reserve constraints are defined by (11) and (12). The renewable energy curtailment is described in (13). Considering the nature of the proposed problem and low voltage characteristics of microgrid, the line limit and network loss are neglected in this paper [31-33].

$$o_{t,l}P_l^{con,\min} \leq P_{t,l}^{con} \leq o_{t,l}P_l^{con,\max} \quad (\forall l \in L, \forall t \in T) \qquad (4)$$

$$\begin{aligned} &-o_{t-1,l} + o_{t,l} - o_{k,l} \leq 0 \\ &(1 \leq k-(t-1) \leq MU_l, \forall l \in L, \forall t \in T) \end{aligned} \qquad (5)$$

$$\begin{aligned} &o_{t-1,l} - o_{t,l} + o_{k,l} \leq 1 \\ &(1 \leq k-(t-1) \leq MD_l, \forall l \in L, \forall t \in T) \end{aligned} \qquad (6)$$

$$o_{t-1,l} - o_{t,l} - v_{t,l} \leq 0 \quad (\forall l \in L, \forall t \in T) \qquad (7)$$

$$-o_{t-1,l} + o_{t,l} - u_{t,l} \leq 0 \quad (\forall l \in L, \forall t \in T) \qquad (8)$$

$$\begin{aligned} P_{t,l}^{con} - P_{t-1,l}^{con} &\leq (2 - o_{t-1,l} - o_{t,l})P_l^{con,\min} \\ &+ (1 + o_{t-1,l} - o_{t,l})RU_l \quad (\forall l \in L, \forall t \in T) \end{aligned} \qquad (9)$$

$$\begin{aligned} P_{t-1,l}^{con} - P_{t,l}^{con} &\leq (2 - o_{t-1,l} - o_{t,l})P_l^{con,\min} \\ &+ (1 - o_{t-1,l} + o_{t,l})RD_l \quad (\forall l \in L, \forall t \in T) \end{aligned} \qquad (10)$$

$$\sum_{l=1}^{L} P_{t,l}^{con} + P_t^{wind} + P_t^{solar} = D_t \quad (\forall t \in T) \qquad (11)$$

$$E_t^{\xi} = P_{w,t}^{wind,\xi} + P_{s,t}^{solar,\xi} - P_t^{wind} - P_t^{solar} \quad (\forall t \in T) \qquad (12)$$

$$\sum_{l \in L} o_{t,l}P_l^{con,\max} \geq R_t + D_t \quad (\forall t \in T) \qquad (13)$$

### B. Flexible demand resource model

As mentioned in Section III, the flexible demand resource can be divided into two types: power-type flexible demand resource and energy-type flexible demand resource. In this paper, $k_1$ represents the power-type flexible demand resource, while $k_2$ represents the energy-type flexible demand resource.

The mathematical model of power-type one is formulated as shown in the follows[25]:

$$D_{K_1,t} = \sum_{k_1=1}^{K_1} f_{k_1,t}\alpha_{k_1,t} \quad (\forall t \in T_{k_1}^{op}) \qquad (14)$$

The following constraints are incorporated in this paper to model the energy-type flexible demand resources. The capacity constraint of flexible demand resource is defined in (15). Constraints (16) and (17) represent the ramping-up/down constraints of flexible demand resource. Constraints (18) and (19) represent the min/max continuous on/off time of flexible demand resources.

$$\begin{aligned} \alpha_{k_2,t}D_{k_2}^{\min} &\leq D_{k_2,t} \leq \alpha_{k_2,t}D_{k_2}^{\max} \\ &(\forall k_2 \in K_2, \forall t \in T_{k_2}^{op}) \end{aligned} \qquad (15)$$

$$\begin{aligned} D_{k_2,t} - D_{k_2,t-1} &\leq RU_{k_2} \\ &(\forall k_2 \in K_2, \forall t \in T_{k_2}^{op}) \end{aligned} \qquad (16)$$

$$D_{k_2,t-1} - D_{k_2,t} \le RD_{k_2} \quad (\forall k_2 \in K_2, \forall t \in T_{k_2}^{op}) \tag{17}$$

$$\sum_{\tau=t-T_{k_c}^{\min}}^{t-1} \alpha_{k_2,\tau} - (\alpha_{k_2,t-1} - \alpha_{k_2,t})T_{k_c}^{\min} \ge 0 \quad (\forall k_2 \in K_2, \forall t \in T_{k_2}^{op}) \tag{18}$$

$$\sum_{\tau=t-T_{k_c}^{\max}-1}^{t-1} (1-\alpha_{k_2,\tau}) \le T_{k_c}^{\max} \quad (\forall k_2 \in K_2, \forall t \in T_{k_2}^{op}) \tag{19}$$

Moreover, the "energy-neutral" constraint for both types of demand resource are defined as following:

$$\sum_{t=1}^{T} D_{k_c,t} = \sum_{t=1}^{T} f_{k_c,t} \quad (\forall k_c \in K_c) \tag{20}$$

*C. Economic benefit measurement and risk mitigation model*

As there are numbers of uncertain factors involved in the microgrid operation, the realized economic benefit induced by flexible demand resources may vary under different scenarios. Over-compensation to the demand resource providers may cause unnecessary financial loss for the microgrid operator, which could eventually jeopardize the sustainability of incentive-based DSM programs. To address this issue, a chance constraint is proposed in this paper to manage the risks of over-compensation. As shown in (21)-(22), the probability that economic benefit is greater than compensation should be larger than $1-\varepsilon$ (confidence level). Consequently, the per unit compensation $\pi_{kc}$ should be equal to the total compensation (i.e. economic benefit) divided by the corresponding controllable capacity of demand resource which indicates its participating capacity in DSM programs. In this way, the chance of financial loss for the microgrid operator is bounded.

$$Pr[(Benefit - \pi_{k_c} \cdot Capacity_{k_c}) > 0] \ge 1-\varepsilon \ (\forall k_c \in K_c) \tag{21}$$

Where

$$Benefit = \min Cost^{withoutDSM} - \min Cost^{withDSM} \tag{22}$$

*D. Iteration-based chance-constrained method for flexible demand resource pricing*

As indicated by (22), the economic benefit is measured by the operation cost deviation between adopting or not adopting the incentive-based DSM program. It can be derived by calculating the deviation between two optimal values, which is quite complicated to be solved directly since it is an optimization problem with two minimization sub-problems. To effectively calculate the economic benefit and determine the compensation, an iteration-based chance-constrained method is proposed in this paper, as shown in the following.

**Algorithm 1 Iteration-based chance-constrained method for flexible demand resource pricing**

1: **Initialize:** Generate N scenarios for RES
2: **for** $n = 1,2…N$ **do**
3:     **Minimize** *Cost1* without DSM
4:     **Minimize** *Cost2* with DSM
5:     **Calculate** *Benefit(n)* = *Cost1* - *Cost2* of scenario *n*
6: **end for**
7: **Sort** the value of each *Benefit(n)*
8: **Pick** the [$(1-\varepsilon)\times N$]th *Benefit* to determine the compensation $\pi$

$$\pi = \frac{Benefit[(1-\varepsilon)\times N]}{Capacity_{[(1-\varepsilon)\times N]}}$$

With the algorithm above, the original chance-constrained optimization problem has been transformed into N deterministic sub-questions with different scenarios. Each sub-question solves the microgrid optimal operation model twice, as shown in step 3 and 4. Each time the problem can be modeled as a mixed integer linear programing problem, which can be efficiently solved by commercial off-the-shelf solvers. To meet certain confidence level, the algorithm picks the [$(1-\varepsilon)\times N$]th scenario after sorting by the value of economic benefit to determine the compensation distributed to the flexible demand resource providers.

It should be mentioned that the economic benefit derived based on the chance-constrained method is completely allocated to the flexible demand resource providers for compensation. The idea behind this design is that the microgrid operator is willing to encourage participation of incentive-based DSM programs as much as possible to improve the overall efficiency. Meanwhile, the microgrid operator does not want to cause financial losses due to over-compensation considering the uncertainty risks involved in this process. The proposed iteration-based chance-constrained method could help the operator achieve those two targets simultaneously. Moreover, the economic benefit shared with resource providers can be dynamically adjusted by selecting different chance-constrained criteria, correspondingly the trade-off between those two goals can be managed.

## V. CASE STUDY

In this section, the proposed flexible demand resource pricing scheme is illustrated and examined in a sample microgrid. The models are coded in MATLAB and solved by the solver CPLEX 12.7.1. All the experiments are implemented on a computer with Intel(R) Core (TM) i5-8300H CPU@2.30 GHz and 8 GB memory.

*A. Simulation Setup*

Without loss of generality, a sample microgrid with three conventional generators, one wind farm and one solar station, is assumed as the testbed. The single-line diagram of microgrid system is shown in Fig. 4 according to [32]. The parameter settings of conventional generators are listed in Table I and II according to [31]. The output data of renewable generation is considered as uncertain factors in each scenario, generated by the Monte Carlo simulation [29]. The wind power and the solar power are assumed to range according to the output of renewable generation of the typical day in summer of ERCOT. In this case, the number of scenarios is set to be 100.

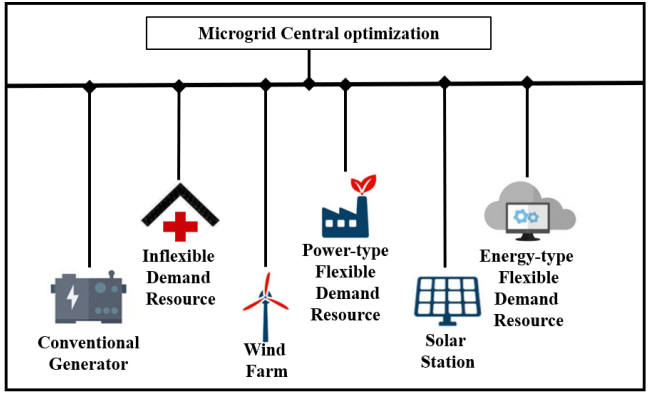

Fig. 4. Architecture of sample microgrid

TABLE I
Conventional Generators Settings

| Unit | Max/Min Output (kW) | Ramp Up/Down Rate (kW/h) | Min Up/Down Time (h) |
|---|---|---|---|
| G1 | 120/30 | 60/60 | 3/3 |
| G2 | 200/80 | 120/120 | 1/1 |
| G3 | 80/15 | 35/35 | 2/2 |

TABLE II
Fuel Settings

| Unit | $a$ ($) | $b$ ($/kWh) | $c$ ($/kW$^2$h) | Start-up Cost ($) |
|---|---|---|---|---|
| G1 | 1.5556 | 0.1489 | 0.0016 | 5.0 |
| G2 | 1.6800 | 0.1667 | 0.0029 | 4.5 |
| G3 | 0.3444 | 0.1813 | 0.0013 | 3.0 |

As for the demand side, the data of demand resource in the typical day in summer is obtained from a report of Everbright Securities [34]. As mentioned in Section III, flexible demand resource can be divided into six categories. Considering the fact that flexible demand resources is only a small percentage of resources in most cases[35], the proportion of flexible demand resource is assumed to be 10% of total installed capacity. The detailed parameter settings for flexible demand resources are illustrated by Fig. 5. and Table III. The parameter of power-type flexible demand resource is partially captured from Fig. 3.

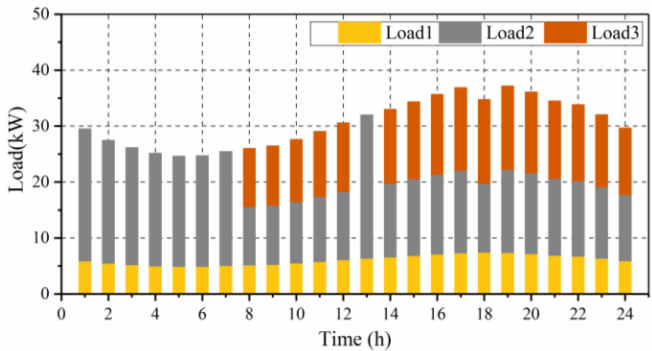

Fig. 5. Profiles of flexible demand resources

TABLE III
Flexible Demand Resource Settings

| Load Group | Max/Min Power (kW) | Ramp Up/Down Rate (kW/h) | Min/Max Continues Time (h) |
|---|---|---|---|
| L1 | 12/5 | 6/6 | 3/- |
| L2 | 45/10 | 20/20 | 4/- |
| L3 | 25/0 | 15/15 | 3/15 |

Besides, other parameters involved in the pricing scheme are reported as follows. The spinning reserve of the microgrid is set to be 25kW. The penalty for curtailed wind and solar power is set to be 0.6$/kWh. To mitigate the financial risk of over-compensation for demand resource providers, 100 scenarios are considered and the confidence level for chance constraint is set as 85%.

B. *Numerical Results*

1) *Energy-Type DSM*

As the economic benefit of flexible demand resource is calculated based on the deviation between cost with DSM and cost without DSM, the optimal operation results for both cases are reported in Fig. 6-8. As mentioned above, the chance-constrained problem has been transformed into N deterministic sub-questions in different scenarios. Fig. 6-8 show the result of [(1-$\varepsilon$)×N]th scenario after sorting by the value of benefits, which determines the compensation distributed to the customers (i.e. the economic benefit shared with demand resource providers).

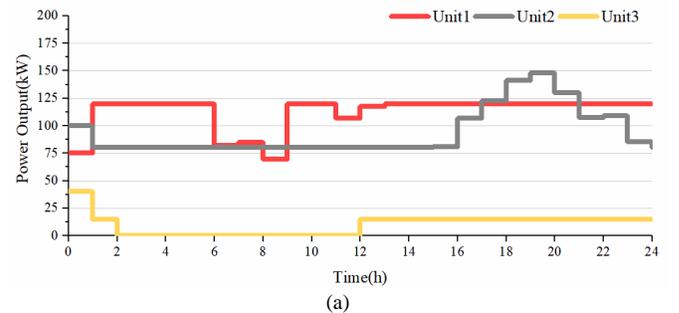

(a)

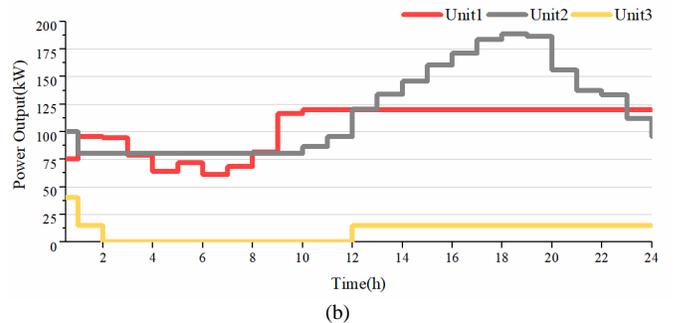

(b)

Fig. 6. (a) Power output of conventional generators with DSM
(b) Power output of conventional generators without DSM

Fig. 6. (a) and Fig. 6. (b) show the scheduled outputs of conventional generators in the cases with DSM and without DSM. Comparing (a) and (b), it is obvious that in valley hours (00:00-9:00), scheduled output without DSM is less than that with DSM while opposite results can be observed in the valley hours (10:00-20:00). This is because that the flexible demand resources reduce the demand in peak hour by shifting loads into valley period to achieve a more economical operation schedule.

Fig. 7. reports the scheduled decisions for energy-type flexible demand resource along with scheduled power output of generation resources. It can be observed that the controlling signal for flexible demand resources (i.e. Δload) in the peak hour (14:00-20:00) is negative. To explain the decisions more clearly, RES outputs, original and adjusted load curves of flexible demand resources are compared in Fig. 8. The trend of grey curve which represents the adjusted load curve for flexible demand resources is similar to blue and green column representing RES outputs. Considering the low variable cost characteristic of RES, it becomes more economical for the microgrid operator to schedule more demands in the hours with more RES outputs, as shown in Fig. 8.

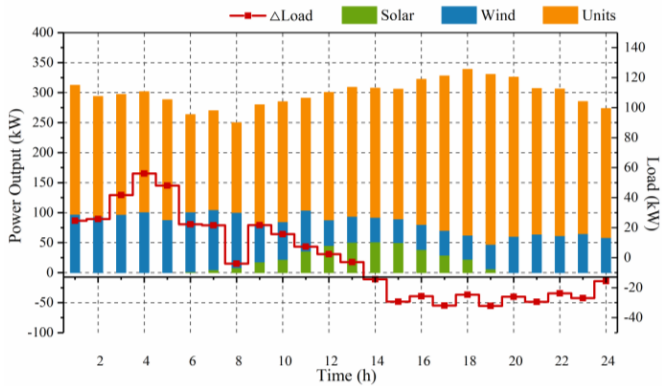

Fig. 7. Scheduled decisions for energy-type flexible demand resource

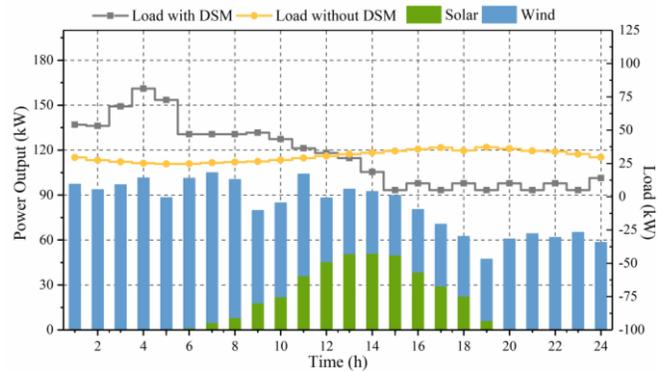

Fig. 8. Results comparison between original and adjusted load curves for flexible demand resources

Apart from flexible demand resource available in the entire operation time span, different types of flexible demand resource are discussed as follows. Fig. 9-10 illustrate scheduled decisions for energy-type flexible demand resources only available for the peak period or valley period. The controlling signal for flexible demand resources (i.e. Δload) is indicated by the red curves in those figures.

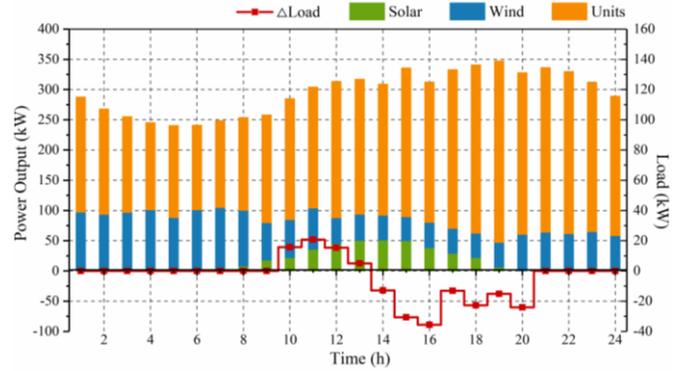

Fig. 9. Scheduled decisions for energy-type flexible demand resource available for peak period

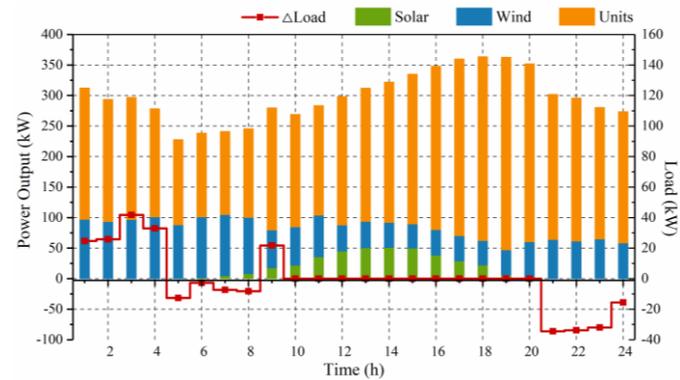

Fig. 10. Scheduled decisions for energy-type flexible demand resource available for valley period

*2) Power-Type DSM*

The characteristics of power-type flexible demand resource are assumed partially based on Fig. 3. The scheduled decisions for the power-type resources with different availability time are reported in Fig. 11-13. It is obvious that the specific power curve of power type flexible demand resource is shifted to the hours with more RES outputs. For instance, power-type flexible demand resource available in full time is shifted from the peak hour to the valley hour. Similar results of flexible demand resource available for peak period and for valley period are reported in Fig. 12-13.

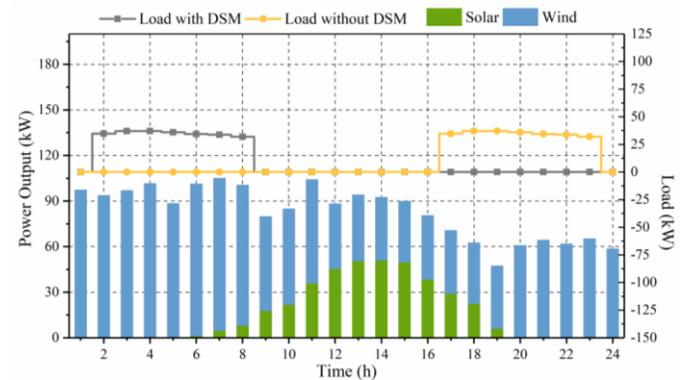

Fig. 11. Scheduled results of power-type flexible demand resource available in full time

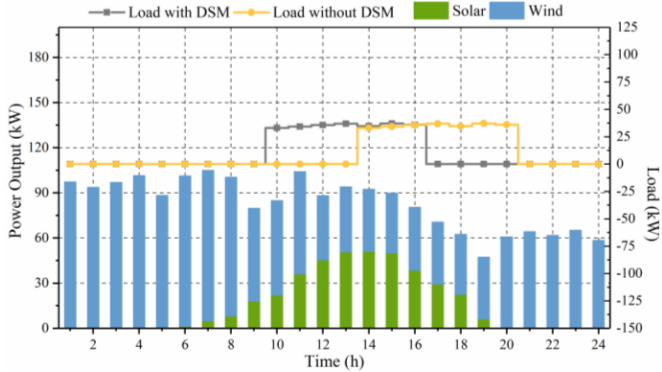

Fig. 12. Scheduled results of power-type flexible demand resource available for peak period

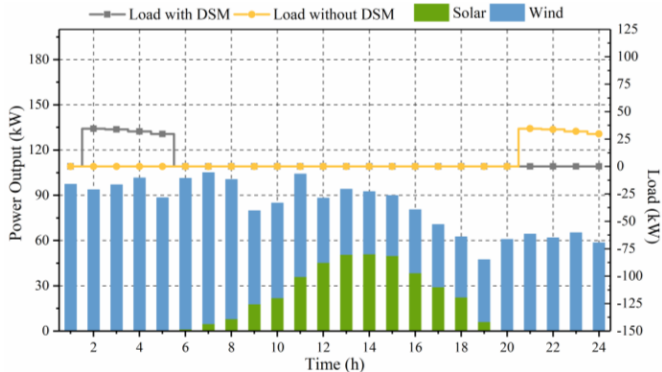

Fig. 13. Schedule results of power-type flexible demand resource available for valley period

*3) Pricing Scheme*

Based on the iteration-based chance-constrained method, the economic benefit and corresponding compensation for different types of flexible demand resources are calculated, as shown in Table IV. All those values in Table IV are calculated for one-day period with a typical operation data. However, the proposed pricing scheme can be easily extended to longer period, such one week or one month, if corresponding data is available. The variation of operation data is captured by multiple scenarios. For each scenario, the benefit is calculated by the deviation between operation cost without DSM and that with DSM, based on (1)-(20). The compensation value is determined based on (21) and (22), guaranteeing that the probability that economic benefit is greater than compensation should be larger than 1-ε.

According to Table IV, it can be observed that energy-type flexible demand resource available in the entire operation time span generate the most economic benefit. This is intuitive since it has the largest flexibility. Consequently, demand resource providers of this type of flexible demand resource receive most compensations with 1.71$/kW per day. In contrast with it, per unit compensation for demand resource available in peak period and valley period are lower. Moreover, compensation for resources available in the valley period is greater than that in the peak period. This can be explained that demand resources contribute greater value during valley hours. Meanwhile, the compensation for power-type flexible demand resource follows the similar pattern, while the overall level of compensation for power-type flexible demand resource is lower than that of energy-type flexible demand resource.

The computational time for each type of flexible demand resource is calculated in Table IV, either. With 100 scenarios, the computational time is less than 30 seconds, which is sufficiently fast for a pricing scheme.

TABLE IV
Compensations for Flexible Demand Resource

| Type | Available Time | Benefits ($/d) | Compensation ($/kW·d) | Computational time (sec) |
|---|---|---|---|---|
| energy-type | Full time | 68.46 | 1.711 | 28.22 |
| energy-type | Peak hour | 64.27 | 1.612 | 19.30 |
| energy-type | Valley hour | 65.49 | 1.690 | 20.80 |
| power-type | Full time | 16.97 | 0.424 | 18.02 |
| power-type | Peak hour | 6.68 | 0.167 | 16.25 |
| power-type | Valley hour | 8.04 | 0.212 | 16.12 |

The economic benefit sharing between end-users and the microgrid operator can be dynamically adjusted with different iteration-based chance-constrained criterion, as shown in Fig. 14. It can be observed that the economic benefit allocated to the demand resource providers decrease as the chance-constrained criterion increases. That is because the chance-constrained criterion represents the probability of avoiding financial loss for the microgrid operator. Therefore, a more risk-averse microgrid operator can reduce the risk of financial loss by flexibly adjusting the chance-constrained criterion. Corresponding, the economic benefit shared with demand resource providers will be reduced.

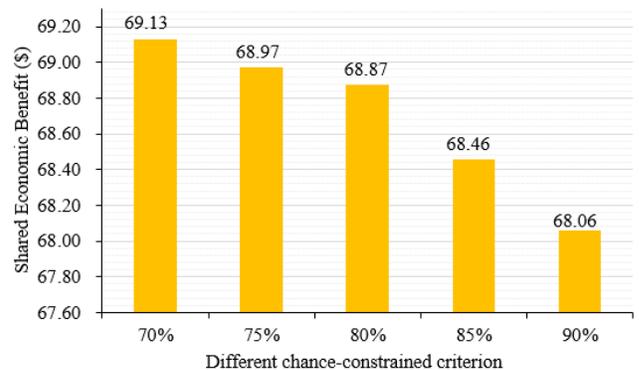

Fig. 14. Economic benefit shared with flexible demand resources at different risk criterion

*4) Scalability Analysis*

In the previous cases, the flexible demand resources compose 10% of total installed capacity to simulate practical cases. To further analyze the impact of flexible demand resource proportion on economic benefit and compensation, the scalability analysis is conducted in this sub-section. As shown in Fig. 15, the economic benefits of flexible demand resources at different scales are reported. It can be observed that the benefits increase as the scale of flexible demand resource increases. However, the benefits do not increase linearly. While the proportion rises to a certain extent, the value of benefits tends to saturation. This can be explained that most of the economic benefit comes from the improvement for utilization of RES. The performance enhancement induced by flexible demand resources would be saturated once there is no improvement space on RES utilization. In this case, as the RES takes 20% of total installed capacity, the economic benefit of flexible demand resource would reach saturation once it increases to 20%.

To further illustrate the analysis above, the scheduled decisions for flexible demand resource and RES output in the case of 20% total installed capacity under the confidence level of 85% are reported in Fig. 16. It can be observed that the controlling signal, Δload is following the pattern of RES output curve to maximize the economic performance.

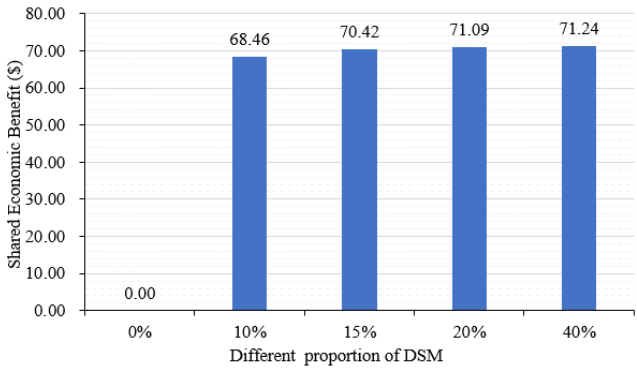

Fig. 15. Economic benefit with flexible demand resources at different scales

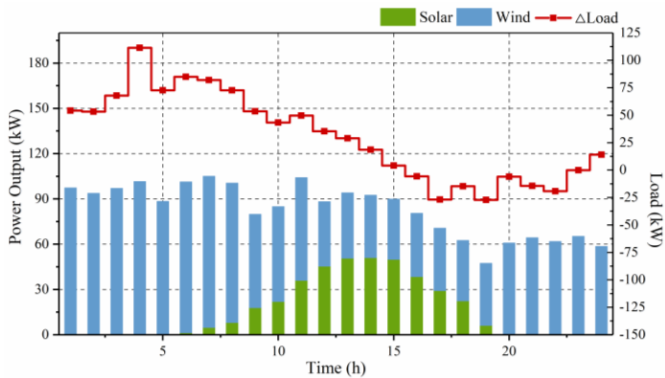

Fig. 16. Scheduled decisions for flexible demand resource and RES output of 20% total installed capacity

The implication of this scalability analysis on the pricing scheme is that the per unit value of flexible demand resources would decrease as its relative scale increasing to a certain level. Therefore, the microgrid operator should dynamically adjust the per unit compensation once the relative scale of flexible demand resources changed.

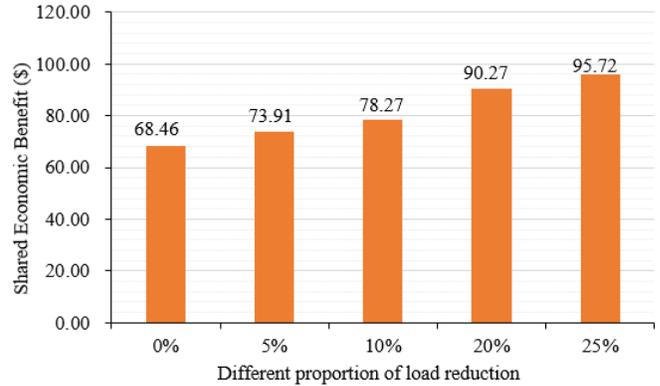

Fig. 17. Economic benefit with different consumption reduction ratio

To further explore the impact of the "energy-neutral" constraint, the economic benefits of flexible demand resources are calculated with different level of consumption reduction during the demand side management process, as shown in Fig. 17. It can be noted that the benefit increases as the consumption reduction increases. This can be explained that the consumption reduction would cause the decrease of power output of units and the operation cost of microgrid operator, making the economic benefit of demand resource providers increase.

VI. CONCLUSION

In this paper, a flexible demand resource pricing scheme is proposed to determine the compensation for incentive-based DSM program without knowing the utility functions of end-users. An iteration-based chance-constrained method is provided to effectively calculate the compensation for demand resource providers and mitigating the financial risks for the microgrid operator. The economic benefit could be dynamically shared between flexible demand resource providers and the microgrid operator by adjusting the chance-constrained criterion. Numerical case studies results demonstrate the effectiveness of the proposed pricing scheme. Also, a scalability analysis is conducted and the results show that the economic benefit mainly comes from performance enhancement of generation due to "energy-neutral" setting and will saturate as the scale of flexible demand resource increases.

REFERENCES

[1] K. Lai, M. S. Illindala, and M. A. Haj-ahmed, "Comprehensive protection strategy for an islanded microgrid using intelligent relays." pp. 1-11.
[2] Q. Qdr, "Benefits of demand response in electricity markets and recommendations for achieving them," *US Dept. Energy, Washington, DC, USA, Tech. Rep*, 2006.
[3] M. Muratori, and G. Rizzoni, "Residential demand response: Dynamic energy management and time-varying electricity pricing," *IEEE Transactions on Power systems,* vol. 31, no. 2, pp. 1108-1117, 2016.
[4] G. Cardoso, M. Stadler, M. C. Bozchalui, R. Sharma, C. Marnay, A. Barbosa-Póvoa, and P. Ferrão, "Optimal investment and scheduling of

- [4] distributed energy resources with uncertainty in electric vehicle driving schedules," *Energy,* vol. 64, pp. 17-30, 2014.
- [5] F. Wang, L. Zhou, H. Ren, X. Liu, S. Talari, M. Shafie-khah, and J. P. Catalão, "Multi-Objective Optimization Model of Source–Load–Storage Synergetic Dispatch for a Building Energy Management System Based on TOU Price Demand Response," *IEEE Transactions on Industry Applications,* vol. 54, no. 2, pp. 1017-1028, 2018.
- [6] L. Zhao, Z. Yang, and W.-J. Lee, "The Impact of Time-of-Use (TOU) Rate Structure on Consumption Patterns of the Residential Customers," *IEEE Transactions on Industry Applications,* vol. 53, no. 6, pp. 5130-5138, 2017.
- [7] A. Safdarian, M. Fotuhi-Firuzabad, and M. Lehtonen, "A medium-term decision model for DisCos: Forward contracting and TOU pricing," *IEEE Transactions on Power systems,* vol. 30, no. 3, pp. 1143-1154, 2015.
- [8] Z.-f. Tan, L.-w. Ju, H.-h. Li, J.-y. Li, and H.-j. Zhang, "A two-stage scheduling optimization model and solution algorithm for wind power and energy storage system considering uncertainty and demand response," *International Journal of Electrical Power & Energy Systems,* vol. 63, pp. 1057-1069, 2014.
- [9] L. Zhu, Z. Yan, W.-J. Lee, X. Yang, Y. Fu, and W. Cao, "Direct load control in microgrids to enhance the performance of integrated resources planning," *IEEE Transactions on Industry Applications,* vol. 51, no. 5, pp. 3553-3560, 2015.
- [10] C. Zhang, Y. Xu, Z. Y. Dong, and J. Ma, "Robust operation of microgrids via two-stage coordinated energy storage and direct load control," *IEEE Transactions on Power Systems,* vol. 32, no. 4, pp. 2858-2868, 2017.
- [11] G. Gutiérrez-Alcaraz, E. Galván, N. González-Cabrera, and M. Javadi, "Renewable energy resources short-term scheduling and dynamic network reconfiguration for sustainable energy consumption," *Renewable and Sustainable Energy Reviews,* vol. 52, pp. 256-264, 2015.
- [12] B. Zhang, and J. Baillieul, "Control and communication protocols based on packetized direct load control in smart building microgrids," *Proceedings of the IEEE,* vol. 104, no. 4, pp. 837-857, 2016.
- [13] N. I. Nwulu, and X. Xia, "Optimal dispatch for a microgrid incorporating renewables and demand response," *Renewable Energy,* vol. 101, pp. 16-28, 2017.
- [14] J. Aghaei, and M.-I. Alizadeh, "Multi-objective self-scheduling of CHP (combined heat and power)-based microgrids considering demand response programs and ESSs (energy storage systems)," *Energy,* vol. 55, pp. 1044-1054, 2013.
- [15] "2018 Assessment of Demand Response and Advanced Metering," F. E. R. Commission, ed., 2018.
- [16] W. Liang, M. Liu, F. Song, W. Wu, K. Zhou, and A. Jin, "Power system dynamic economic dispatch with controllable air-conditioning load groups." pp. 962-967.
- [17] C. D. Korkas, S. Baldi, I. Michailidis, and E. B. Kosmatopoulos, "Occupancy-based demand response and thermal comfort optimization in microgrids with renewable energy sources and energy storage," *Applied Energy,* vol. 163, pp. 93-104, 2016.
- [18] O. Erdinc, A. Taşcıkaraoğlu, N. G. Paterakis, Y. Eren, and J. P. Catalão, "End-user comfort oriented day-ahead planning for responsive residential HVAC demand aggregation considering weather forecasts," *IEEE Transactions on Smart Grid,* vol. 8, no. 1, pp. 362-372, 2017.
- [19] M. H. Imani, P. Niknejad, and M. Barzegaran, "The impact of customers' participation level and various incentive values on implementing emergency demand response program in microgrid operation," *International Journal of Electrical Power & Energy Systems,* vol. 96, pp. 114-125, 2018.
- [20] A. Rabiee, M. Sadeghi, J. Aghaeic, and A. Heidari, "Optimal operation of microgrids through simultaneous scheduling of electrical vehicles and responsive loads considering wind and PV units uncertainties," *Renewable and Sustainable Energy Reviews,* vol. 57, pp. 721-739, 2016.
- [21] W. S. Ho, S. Macchietto, J. S. Lim, H. Hashim, Z. A. Muis, and W. H. Liu, "Optimal scheduling of energy storage for renewable energy distributed energy generation system," *Renewable and Sustainable Energy Reviews,* vol. 58, pp. 1100-1107, 2016.
- [22] B. Zhang, Q. Li, L. Wang, and W. Feng, "Robust optimization for energy transactions in multi-microgrids under uncertainty," *Applied Energy,* vol. 217, pp. 346-360, 2018.
- [23] Y. Li, T. Zhao, P. Wang, H. B. Gooi, L. Wu, Y. Liu, and J. Ye, "Optimal operation of multimicrogrids via cooperative energy and reserve scheduling," *IEEE Transactions on Industrial Informatics,* vol. 14, no. 8, pp. 3459-3468, 2018.
- [24] R. H. Lasseter, "Microgrids." pp. 305-308.
- [25] Q. Li, N. Song, J. Wang, and H. Zhong, "A Pattern and Method of Optimized Power Utilization Based on Consumers' Interaction Capability," *Power System Technology,* vol. 40, no. 6, pp. 1816G1824, 2016.
- [26] M. Ghamkhari, and H. Mohsenian-Rad, "Energy and performance management of green data centers: A profit maximization approach," *IEEE Transactions on Smart Grid,* vol. 4, no. 2, pp. 1017-1025, 2013.
- [27] N. D. a. R. Commission, *Notice on the implementation of Coal-Electricity Price linkage in north China power grid*, 2005.
- [28] A. Zakariazadeh, S. Jadid, and P. Siano, "Stochastic multi-objective operational planning of smart distribution systems considering demand response programs," *Electric Power Systems Research,* vol. 111, pp. 156-168, 2014.
- [29] Z. Ding, and W.-J. Lee, "A stochastic microgrid operation scheme to balance between system reliability and greenhouse gas emission." pp. 1-9.
- [30] Q. Wang, Y. Guan, and J. Wang, "A chance-constrained two-stage stochastic program for unit commitment with uncertain wind power output," *IEEE Transactions on Power Systems,* vol. 27, no. 1, pp. 206-215, 2012.
- [31] J. Zeng, Q. Wang, J. Liu, J. Chen, and H. Chen, "A Potential Game Approach to Distributed Operational Optimization for Microgrid Energy Management With Renewable Energy and Demand Response," *IEEE Transactions on Industrial Electronics,* vol. 66, no. 6, pp. 4479-4489, 2018.
- [32] M. Ross, C. Abbey, F. Bouffard, and G. Jos, " Multiobjective optimization dispatch for microgrids with a high penetration of renewable generation," *IEEE Transactions on Sustainable Energy,* vol. 6, no. 4, pp. 1306-1314, 2015.
- [33] P. Li, D. Xu, Z. Zhou, W.-J. Lee, and B. Zhao, "Stochastic optimal operation of microgrid based on chaotic binary particle swarm optimization," *IEEE Transactions on Smart Grid,* vol. 7, no. 1, pp. 66-73, 2015.
- [34] http://www.hlfdw.com/a/xingyekuaixun/guonaxinwen/20170818/173999.html.
- [35] V. Durvasulu, and T. Hansen, "Benefits of a Demand Response Exchange Participating in Existing Bulk-Power Markets," *Energies,* vol. 11, no. 12, pp. 3361, 2018.